\documentclass{amsart}
\usepackage{amssymb,amsrefs}
\usepackage{graphics}

\theoremstyle{plain}
\newtheorem{lemma}{Lemma}[section]
\newtheorem{theorem}[lemma]{Theorem}
\newtheorem{corollary}[lemma]{Corollary}

\newcommand{\Action}{\mathcal{A}}
\newcommand{\quotient}{/}
\newcommand{\Xt}{X_{H_t}}

\newcommand{\R}{\mathbb{R}}
\newcommand{\C}{\mathbb{C}}

\newcommand{\Z}{\mathbb{Z}}

\newcommand{\0}{{x_{0}}}

\newcommand{\infinity}{\infty}

\newcommand{\die}[2][]{\frac{\partial {#1}}{\partial {#2}}}

\newcommand{\pf}{{}_{*}}
\newcommand{\pb}{{}^{*}}
\renewcommand{\d}{{\,d}}

\newcommand{\e}{\operatorname{e}}
\newcommand{\Exp}[1]{{\e^{#1}}}

\newcommand{\FlowProj}{P}
\newcommand{\Sphi}{\tilde \phi}
\newcommand{\iext}{\tilde \imath}

\begin{document}
\begin{abstract}
This paper introduces techniques of symplectic topology to the study of homoclinic orbits in Hamiltonian systems. The main result is a strong generalization of homoclinic existence results due to S\'er\'e and to Coti-Zelati, Ekeland and S\'er\'e \cites{Sere95,CZES}, which were obtained by variational methods. Our existence result uses a modification of a construction due to Mohnke 
\cite{Mohnke01} (originally in the context of Legendrian chords), and an energy--capacity inequality of Chekanov \cite{Chekanov98}.  In essence, we show the existence of a homoclinic orbit by showing a certain Lagrangian embedding cannot exist.

We consider a (possibly time dependent) Hamiltonian system on an exact symplectic manifold
$(W, \omega = d \lambda)$ with a hyperbolic rest point.  

In the case of periodic time dependence, we show the existence of an orbit homoclinic to the rest point if $\lambda(X_H) - H $ is positive and proper, $H$ is positive outside a compact set and proper, and $(W, \omega)$ admits the structure of a Weinstein domain.

In the autonomous case, we establish the existence of an orbit homoclinic to the rest point if the critical level is of restricted contact-type, and the critical level has a Hamiltonian displaceable neighbourhood.
\end{abstract}
\title{Homoclinic orbits and Lagrangian embeddings}
\author[S. Lisi]{Samuel T. Lisi}
\address{Universit\'e de Montr\'eal, Montr\'eal (Qu\'ebec), Canada}
\email{lisi@dms.umontreal.ca}
\date{11 August 2007}

\maketitle

\section{Introduction}

In this work, we establish the existence of homoclinic orbits for certain classes of  Hamiltonian systems.  Similar results in this direction were obtained by Coti-Zelati, Ekeland and S\'er\'e \cite{CZES} for time-dependent systems (with periodic time dependence), and by Hofer and Wysocki \cite{HWhomo90} and by S\'er\'e \cite{Sere95} in the autonomous case.  Their results use variational methods and global analysis.  In this work, the result is an application of methods from symplectic topology, notably of an energy--capacity inequality of Chekanov \cite{Chekanov98}, the basic idea of which goes back to work of Floer and of Gromov with pseudoholomorphic disks in symplectic manifolds \cites{FloerLagrangian,Gromov85}.

Given a symplectic manifold $(M, \omega)$ and a $C^2$ function $H : S^1 \times M \to \R$, we define the (time dependent) \textit{Hamiltonian vector field} $\Xt$ by $i_{\Xt} \omega = - d H_t$, for each $t \in S^1 = \R \quotient \Z$.   In the case that $H$ is independent of $t$, we say it is an autonomous (or time independent) Hamiltonian.  We let $\phi_t$ be the flow of this dynamical system.

If the function $H(t, \cdot)$ has a critical point at $\0 \in M$ for all $t \in S^1$, 
the point $\0$ is a fixed point of the flow.  
Any solution which is asymptotic to $\0$ in both forward and backward time 
is said to be an orbit \textit{homoclinic} to $\0$.  

The question of the existence of homoclinic orbits traces back to work of Poincar\'e.  
Most recent results have used variational methods (as, for instance, in the articles by S\'er\'e and by Coti-Zelati, Ekeland and S\'er\'e, and also as discussed in Rabinowitz's survey article \cite{RabinowitzSurvey}). 

In this article, we will use symplectic topology methods.  These methods have proved surprisingly effective at answering existence questions for periodic orbits, starting with Floer's work on the Arnol{\cprime}d Conjecture \cite{FloerArnoldConjecture}.  The key point of the argument will be to adapt to the homoclinic setting a construction originally used by Mohnke in the context of Reeb Chords (``brake'' orbits) \cite{Mohnke01}.  This permits us to transform the question of existence of homoclinic orbits into one of existence of Lagrangian embeddings.  This latter question is well studied in symplectic topology by means of pseudoholomorphic disks.  

In order to state the results, we will review a number definitions from symplectic topology.  

A symplectic manifold $(W, \omega)$ is \textit{exact} if $\omega = d \lambda$, for a one form $\lambda$ on $W$.  We say $(W, \omega, Y, F)$ is a \textit{Weinstein domain} if $(W, \omega)$ is a symplectic manifold, $L_Y \omega = \omega$, $F$ is an exhausting Morse function on $W$, and $Y$ is gradient-like for $F$.  Furthermore, if $W$ is non-compact, we require the vector field $Y$ to be complete.  If $W$ is compact, we require the boundary to be a regular level of $F$.  In particular, a Weinstein domain is an exact symplectic manifold, with a primitive given by $i_Y \omega$. 

Following Gromov \cite{Gromov85}, we say that the non-compact symplectic manifold $(W, \omega)$ is geometrically bounded if there exists an almost complex structure $J$ such that  $\omega( \cdot, J \cdot )$ is a Riemannian metric on $W$ whose sectional curvature is bounded above and whose injectivity radius is bounded below.  We also refer to such almost complex structures as having bounded geometry.

A diffeomorphism $\Phi : W \to W$ is said to be a 
Hamiltonian diffeomorphism if it is the time $1$ map of a time dependent Hamiltonian flow.
To such a $\Phi$, we may associate the Hofer norm, denoted $|| \Phi ||$.  This is given by taking the infimum of the oscillation over all Hamiltonians that generate $\Phi$ :
\[
|| \Phi || = \inf_{H_t \text{ generates } \Phi} \, \int_{0}^{1} (\sup_{W} H_t - \inf_{W} H_t ) \d t.
\]

A set $S \subset W$ is said to be Hamiltonian displaceable if there exists a (time dependent) 
Hamiltonian $F : [0, 1] \times W \to \R$, with compact support, so that the time one map $\Psi$ of the flow of $X_F$ displaces the closure of $S$, i.e. $\overline {\Psi(S) }\cap \overline S = \emptyset$.  

A hypersurface $M$ in an exact symplectic manifold $(W, \omega = d \lambda)$
is of \textit{restricted contact-type} if $\lambda|_M$ is a contact form on $M$.  In the special case of a hypersurface given as a regular level set of a function $F$, 
this is equivalent to having $\lambda(X_F) \ne 0$ on $M$.

By analogy, we say that a Hamiltonian vector field $\Xt$ with periodic time dependence is \textit{Reeb-like} if 
$\lambda(\Xt) - H_t > 0$.  Equivalently, this is the condition that $\Lambda = \lambda - H_t dt$ is a contact form on $S^1 \times W$.

We establish the following two results :
\begin{theorem} \label{T:Autonomous}
Let $(W, \omega = d \lambda)$ be an exact symplectic manifold of dimension $2n$, of bounded geometry, 
and $H : W \to \R$ an autonomous Hamiltonian with $H(\0) = 0$ and $dH(\0) = 0$, and with $\0$ a hyperbolic zero of $X_H$.  

Suppose that $H^{-1}(0) \subset W$ is compact and of \textit{restricted contact--type} away from $\0$, i.e. $\lambda(X_H) > 0$ on $H^{-1}(0) \setminus \{ \0 \}$.

If for any neighbourhood $U$ of $\0$, $H^{-1}(0) \setminus U$ is Hamiltonian displaceable in $W$, 
then there exists an orbit of the Hamiltonian vector field $X_H$ homoclinic to $\0$.
\end{theorem}

\begin{theorem} \label{T:SelfCalibrated}
Let $(M, \omega, Y, F)$ be a Weinstein domain of dimension $2n$, such that the vector field $Y$ is complete.  Let $\lambda = i_Y \omega$.

Let $H : S^1 \times M \to \R$ be a time dependent Hamiltonian, with periodic time dependence ($S^1 = \R \quotient \Z$).  Let $\Xt$ denote the associated time dependent Hamiltonian vector field. 

Suppose $H(t, \0) \equiv 0$, $dH(t, \0) \equiv 0$ and $\Xt$ has a hyperbolic zero at $\0$ for all $t \in S^1$.  
Assume that the Hamiltonian system is Reeb-like away from $\0$, i.e.   
\begin{equation} \label{E:selfcalibrating}
\lambda( \Xt) - H_t > 0 \quad \text{ on } M \setminus \{ \0 \}.
\end{equation}
and that both $H(t, \cdot)$ and $\lambda(\Xt) - H_t$ are proper for each $t \in S^1$ and that $H$ is positive outside a compact set.

Then, there exists an orbit homoclinic to $\0$.
\end{theorem}

In particular, by Theorem \ref{T:Autonomous}, taking $W = \R^{2n}$, we obtain the result due to S\'er\'e \cite{Sere95}.  In the non-autonomous case, specialising to the case $M = \R^{2n}$, Theorem \ref{T:SelfCalibrated} generalizes the class of Hamiltonians for which the existence result of Coti-Zelati, Ekeland and S\'er\'e \cite{CZES} holds, but we do not obtain their multiplicity result.  (We note that their super-quadratic growth condition implies the conditions of the theorem.)  
Cieliebak and S\'er\'e \cite{CieliebakSereMultibump, CieliebakSereHomoclinic} have also generalized \cite{CZES} :
 they give a lower bound for the number of homoclinic orbits 
under conditions somewhat more restrictive than the ones we consider here.  
In particular, they work on a cotangent bundle,  and require a form of super-quadratic growth on $\lambda(\Xt) - H_t$.  In the time dependent case, they show there are infinitely many unparameterized homoclinic orbits.

A simple example of a Hamiltonian system in which Theorem \ref{T:SelfCalibrated} is applicable is a generalization of the superquadratic Hamiltonians considered by \cite{CZES} to the case of a cotangent bundle.
Indeed, let $M$ be a compact manifold.  The cotangent bundle $T\pb M$ has a canonical symplectic structure $\omega = d \Theta = \sum d p_i \wedge d q_i$, where $\Theta = \sum p_i d q_i$ is the Liouville 1--form.  
Let $g$ be a Riemannian metric on $M$, which we then dualize to obtain a metric on $T\pb M$.  
Consider a Hamiltonian $H$ on $T\pb M$ of the following form :
\begin{equation*}
H(t, q, p) = \frac{1}{2} || p ||^2 - V(q) + R(t, q, p),
\end{equation*}
where we take $V : M \to \R$ to be a Morse function with $V \ge 0$ and with a unique minimum $q_0$ at which $V(q_0) = 0$, and $R : S^1 \times T\pb M \to \R$ to be one--periodic, with $R(t, q_0, 0) = 0$, $d R(t, q_0, 0) = 0$ and $d^2 R(t, q_0, 0) = 0$. 
In this case, $(q_0, 0)$ is a critical point of the Hamiltonian, satisfying all the non-degeneracy conditions.
If we take $R$ super-quadratic in $p$, we obtain the case of \cite{CieliebakSereMultibump, CieliebakSereHomoclinic}.  
Consider instead the more general case in which $R$ has super-linear growth : suppose there exists a constant $\mu \ge 1$ so that 
\begin{equation*}
0 \le R(t, q, p) \le \frac{1}{\mu} dR \left [ \sum p_i \frac{\partial}{\partial p_i } \right ].
\end{equation*}
Then, $H$ will  be positive for $|| p ||$ large, and is proper.  Furthermore, 
\begin{equation*}
\begin{split}
\Theta( \Xt ) - H_t &= \frac{1}{2} || p||^2 + V(q) + dR[Y] - R \\
		&\ge \frac{1}{2} || p ||^2 + V(q) + (\mu - 1) R\\
		& \ge \frac{1}{2} || p ||^2 + V(q).
\end{split}
\end{equation*}
Then, $\Theta( \Xt ) - H_t$ is positive away from $(q_0, 0)$ and is proper.
Thus, by Theorem \ref{T:SelfCalibrated},
there exists an orbit homoclinic to $(q_0, 0)$. 

A more elaborate class of examples for which the hypotheses of Theorem \ref{T:Autonomous} are verified comes from Weinstein domains \cite{WeinsteinHokkaido}  : 
\begin{corollary} \label{C:Weinstein}
Let $(W^{2n}, \omega, F, Y)$ be a Weinstein domain.  Suppose there is a unique critical point $\0$ on the level set $F = 0$, of Morse index $n$.  Suppose furthermore that any critical point $y$ with $F(y) < 0$ has Morse index strictly less than $n$.  

If $\0$ is a hyperbolic zero of $X_F$, then there exists an orbit homoclinic to $\0$.
\end{corollary}

In essence, the idea in both the proof of Theorem \ref{T:Autonomous} and in the proof of Theorem \ref{T:SelfCalibrated} is that the unstable manifold of the rest point is a non-compact, immersed, exact Lagrangian submanifold (of a suitably chosen symplectic manifold).  
Using a construction of Mohnke's \cite{Mohnke01}, 
we thicken the unstable manifold in order to construct a closed Lagrangian submanifold (no longer exact).  
In the autonomous case, this thickening is possible by the condition that the level set is of restricted contact type.  In the non-autonomous case, 
the thickening is possible because the Hamiltonian vector field is Reeb-like, 
and because $\lambda(\Xt) - H_t$ is proper.
Then, by Chekanov's result that the displacement energy of a closed Lagrangian is bounded from below by the area of a pseudoholomorphic disk, we will show that non-existence of homoclinic orbits 
prevents a neighbourhood from being Hamiltonian displaceable. 
In the autonomous case, this contradicts our hypothesis.  
In the non-autonomous case, we use the fact that $H_t$ is positive outside a compact set to embed the whole problem in $\C \times M$.  Here, again, we obtain a contradiction from the fact that a compact region in $\C \times M$ can always be displaced by a Hamiltonian diffeomorphism.  

We will first prove Theorem \ref{T:Autonomous}.  In this case already, the ideas behind the construction become apparent.  We will then prove Theorem \ref{T:SelfCalibrated}, emphasizing the main differences from the autonomous case.

\subsection*{Acknowledgements}

I would like to thank Helmut Hofer for many fruitful discussions and for helpful suggestions.  
I would also like to thank Octav Corn\'ea for many helpful suggestions for improving the exposition.
I would finally like to thank the anonymous referee for many outstanding suggestions for improvements, and for meticulous attention to detail.

\section{The autonomous case}

\begin{proof}[Proof of Theorem \ref{T:Autonomous}]

We will argue by contradiction.  Suppose there are no homoclinic orbits on the level $H=0$.  

First, we will show the unstable manifold to the rest point $\0$, $W^u$, is an immersed exact Lagrangian (i.e. $i \pb \lambda = d G$ for some function $G : W^u \to \R$).  Furthermore, we will show that $G$ is proper.

By the hypothesis that $\0$ is a hyperbolic zero of the Hamiltonian vector field $X_H$, the unstable manifold exists.  We recall that the unstable manifold is an immersed $\R^n$.  Let us denote this immersion by $i : \R^n \looparrowright  W$, $i(\R^n) = W^u$, with $i(0) = \0$.  Furthermore, for any compact set $K \subset \R^n$, $i|_K : K \to W$ is an embedding.  
Note that the stable/unstable manifolds of a rest point in a Hamiltonian system are isotropic submanifolds of $W$.
Observe that if $\mu$ is an eigenvalue of the linearization of $X_H$ at $\0$, then $- \mu$ is also.  
Since, by hypothesis, there are no purely imaginary eigenvalues, it follows that both the stable and unstable manifolds, $W^s$ and $W^u$, are $n$-dimensional and hence are Lagrangian.
Thus, $i$ is a Lagrangian immersion : $i\pb \omega =0$.

Observe that $i\pb\lambda$ is a closed form on $\R^n$ since $d( i\pb \lambda) = i \pb d\lambda = i \pb \omega = 0$, and hence is exact.  
Thus, $i\pb \lambda = d G$ for some function $G : \R^n \to \R$.  We may set $G(0) = 0$.
The vector--field $X_H$ is tangent to $W^u$, and so we may pull $X_H|_{W^u}$ back to $\R^n$ using $i^{-1}$.  We will write $Z = i^{-1} \pf X_H|_{W^u}$.  
We have $\d G[ Z ] = \lambda( X_H ) \ge 0$, with equality only at $0$. 
Let $\hat \phi_t$ denote the flow of $Z$ on $\R^n$.  Then, $i \circ \hat \phi_t = \phi_t \circ i$.
We observe furthermore that for all $t \in \R$ :
\begin{equation} \label{E:action}
G( \hat \phi_t(x) ) - G(x) = \int_0^t \frac{d}{d\tau} G(\hat \phi_\tau(x)) d\tau  
		=  \int_0^t \lambda( X_H )|_{\phi_\tau (x) } d \tau.
\end{equation}
Thus, $G$ has no critical points away from $0$, and level sets of $G$ are transverse to $Z$.  Thus, $G(x) \ge 0$ with $G = 0$ only at $0$.

First, we have that $G^{-1}(c)$ is a sphere for small, positive values of $c$.  Indeed, 
by the Stable/Unstable manifold theorem, there is a positive quadratic form $N(x)$ on $\R^n$ so that in a neighbourhood of $0$, the norm $N(x)$ is increasing along orbits of $Z$.  Consider the sphere $\mathcal{S} := \{ N(x) = \delta_0 \}$ for $\delta_0$ sufficiently small.  Let $c$ be the minimum of $G$ on this sphere.  Then $c > 0$ and $G^{-1}(c)$ is contained in $\{ x \, | \, N(x) \le \delta_0 \}$ and thus is compact.  Furthermore, $G^{-1}(c)$ is diffeomorphic to $\mathcal{S}$, where the diffeomorphism is given by following the flow lines of $Z$.  
It follows that $G^{-1}(c)$ is a sphere for $c>0$ sufficiently small.  

We will now show that $G$ is proper by using the assumption of non-existence of periodic orbits.
Let $0 < a < b$.    
By assumption, there are no orbits homoclinic to $\0$.  Thus, there exists a neighbourhood $U$ of $\0$, so 
for all $0 \le t \le b+1$, for all $x \notin U$, $\phi_t(x) \notin U$.  Let $c>0$ be sufficiently small that $G^{-1}(c)$ is a sphere.  By shrinking $U$, we may assume that $G^{-1}(c) \cap U = \emptyset$.  

By the compactness of $H^{-1}(0) \setminus U$, there exists a constant $c_U$ so $\lambda(X_H) \ge c_U > 0$ on $H^{-1}(0) \setminus U$.  
Let $\hat U = i^{-1}( U )$.  Then, by Equation \eqref{E:action}, for 
$x \in \R^n \setminus \hat U$ and $t > 0$  :
\begin{equation} \label{E:superlinear}
G( \hat \phi_t(x) ) - G(x) \ge c_U t.
\end{equation}
Observe also that for any $x \in \R^n$, 
\begin{equation} \label{E:limit}
\lim_{t \rightarrow - \infinity} G(\hat \phi_t(x)) = 0.
\end{equation}

The properness of $G$ now follows from these two facts.
Indeed, consider a sequence of points $x_n \in \R^n$ so that $G(x_n)$ is in $[a,b]$ and converges :
\begin{equation*}
\lim_{n \to \infinity} G(x_n) = g \in [a, b].
\end{equation*}
Then, by \eqref{E:limit}, there exists a sequence of times $t_n \in \R$ so that 
$G(\hat \phi_{t_n}(x_n) ) = c$.  
By \eqref{E:superlinear}, the sequence $t_n$ is bounded, independently of $n$.
From above, $G^{-1}(c)$ is a sphere, and hence, compact.  
Thus, there exists a subsequence $n_k$, a real number $t_*$ and a point $y \in G^{-1}(c) \subset \R^n$ 
so that  :
\begin{equation*}
\hat \phi_{t_{n_k}} ( x_{n_k} )  \to y \quad \text{ and } \quad t_{n_k}  \to t_*.
\end{equation*}
Observe that the map $(t, x) \mapsto \hat \phi_t(x)$ is continuous.  Thus, 
\begin{equation*}
\lim_{k \to \infinity} x_{n_k} 
	= \lim_{k \to \infinity} \hat \phi_{-t_{n_k}} \left ( \hat \phi_{t_{n_k}} ( x_{n_k} ) \right ) 
	= \hat \phi_{- t_*} (y).
\end{equation*}
The properness of $G$ now follows.

We will now show that there exists a positive constant $S$ and a bijective immersion $\iext$ of $[-S, 0] \times ( W^u \setminus U ) \looparrowright W$ so that $\iext \pb \lambda = \Exp{s} dG$.
Indeed, since $d \lambda$ is a symplectic form, there exists a vector field $Y$ so that $d\lambda(Y, \cdot) = \lambda$.  By Cartan's formula, $L_Y \lambda = \lambda$.  
Furthermore, $dH [ Y ]  = \lambda(X_H)  > 0$ away from $\0$. 
Hence, $Y$ is transverse to the level set $H^{-1}(0) \setminus \{ \0 \}$.  
Let $\psi_s$ be the flow of the vector field $Y$.  
By the compactness of $H^{-1}(0)$, there exists $S > 0$ so that
\begin{equation*}
\begin{aligned}
\sigma : [-S, 0] \times H^{-1}(0) \setminus U &\longrightarrow W \\
(s, x) &\mapsto \psi_s(x)
\end{aligned}
\end{equation*}
is an embedding.   We now consider the following embedding :
\begin{equation} \label{E:MohnkeEmbed}
\begin{aligned}
\iext : [-S, 0] \times G^{-1}( [a,b] ) &\longrightarrow W \\
(s, x) &\mapsto \psi_s \circ i(x)
\end{aligned}.
\end{equation}  From $L_Y \lambda = \lambda$ and $i\pb \lambda = dG$, it follows that : 
\begin{equation}
\iext\pb \lambda =  \Exp{s} dG.
\end{equation}

On the other hand, by hypothesis, there exists a Hamiltonian diffeomorphism $\Phi$ with compact support and finite Hofer norm $|| \Phi ||$ displacing $H^{-1}(0) \setminus U$.
Thus, there exists $0 < S' \le S$ so that $\Phi$ displaces the image of $[-S', 0] \times (W^u \setminus U)$.  In particular, it displaces the image of $[-S', 0] \times G^{-1}[a, b]$ for all $b > a$.

The result now follows by contradiction to the following result of Mohnke \cite{Mohnke01}, 
a corollary to a deep result of Chekanov \cite{Chekanov98}.
\end{proof}

\begin{theorem}[Mohnke]  \label{T:Mohnke}
Let $(W^{2n}, \omega = d\alpha)$ be an exact symplectic manifold, of bounded geometry.
Let $L$ be a compact manifold of dimension $n$, with boundary.

Suppose $G : L \to \R$ is a Morse function, so that $\partial L = G^{-1}(a) \cup G^{-1}(b)$, where $a$ and $b$ are regular values of $G$.

Suppose there exists an $S > 0$, and an embedding $p : [-S, 0] \times L \to W$ so $p \pb \alpha =  \Exp{s} dG$, where $s$ is the coordinate on $[-S, 0]$.  Suppose furthermore there exists a Hamiltonian diffeomorphism $\Phi : W \to W$, of compact support and with finite Hofer norm $|| \Phi ||$, displacing the image of $p$.  

Then, $(1 - \Exp{-S}) |b - a | \le || \Phi ||$.
\end{theorem}

\begin{proof}[Sketch of proof]

The idea of the proof is to approximate the boundary of the image of $p$ by a smooth Lagrangian.  We will take $a < b$.

For each $\epsilon > 0$, let $\mu_\epsilon : [-S, 0] \times [a, b] \to \R$ so that $0$ is a regular value, and $\mu_\epsilon^{-1}(0)$ is a smooth circle $\mathcal C_\epsilon$ 
approximating the boundary of $[-S, 0] \times [a, b]$,
with $(1 - \Exp{-S})(b - a) - \epsilon \le \int_{\mathcal C_\epsilon} \Exp{s} d t \le (1-\Exp{-S})(b-a)$.
We arrange, furthermore, that $\frac{\partial \mu_\epsilon}{\partial s}(s,t) = 0$ only at $t=a$ or $t=b$.

We now let $\mathcal {\bar L}_\epsilon = \{ (s, x) \, | \, \mu_\epsilon( s, G(x) ) = 0 \}$.  This is a smooth, closed submanifold of $[-S, 0] \times L$.  
Let $\mathcal L_\epsilon = p( \mathcal {\bar L}_\epsilon )$ be the image under the embedding $p$.  
This is a closed Lagrangian submanifold of $(W, \omega)$.  Furthermore, for any disk $u : D \to W$ with $u( \partial D) \subset \mathcal L_\epsilon$, 
by Stokes' Theorem and construction, there is a $\kappa \in \Z$ so :  
\begin{equation} \label{E:StokesEstimate}
\int_D u \pb \omega = \int_{\partial D} u \pb \alpha = \kappa \left ( (1-\Exp{-S})(b-a) - \epsilon \right ).
\end{equation}
(The integer $\kappa$ is the winding of $u(\partial D)$ around the circle $\mathcal C_\epsilon$.)

By hypothesis, the Lagrangian $\mathcal L_\epsilon$ can be displaced from itself by $\Phi$.  
Thus, by an energy--capacity inequality due to Chekanov \cite{Chekanov98}, we have the existence of a non-constant pseudoholomorphic disk $u$ with boundary in $\mathcal L_\epsilon$ and for which the symplectic area is \textit{a priori} bounded : $\int u \pb \omega \le || \Phi ||$.  A non-constant pseudoholomorphic disk has positive area, so we combine this with \eqref{E:StokesEstimate} to obtain :
\begin{equation}
0 <  (1-\Exp{-S})(b-a) - \epsilon  \le \int_D u \pb \omega \le || \Phi ||.
\end{equation}
This holds for any $\epsilon > 0$, and the result now follows.
\end{proof}

\begin{proof}[Proof of Corollary \ref{C:Weinstein}]

It suffices to verify that for any neighbourhood $U$ of $\0$, $F^{-1}(0) \setminus U$ is Hamiltonian displaceable in a suitable symplectic manifold.  We will first recall some results from symplectic topology on Weinstein domains.

We say that a Weinstein domain of dimension $2n$ is \textit{subcritical} if the Morse indices of the critical points of $F$ have index no greater than $n-1$.  We say it is complete if the flow of $Y$ exists for all time (in particular then, $W$ is non-compact and without boundary).  

There is a natural method to construct a complete Weinstein domain, given a compact one with boundary.  
Suppose $W$ is compact, and its boundary is a regular level of $F$.  From this, we construct a \textit{completion}, which we denote by $(\widetilde W, \tilde \omega, \tilde Y, \tilde F)$.  
We set $M = \partial W$.  The completion $\widetilde W$ is obtained by gluing on the symplectization $\R^+ \times M$ to the boundary of $W$ and extending $Y, F$ and $\omega$ over this.  This now gives a complete Weinstein domain.

By Biran and Cieliebak \cite{CieliebakBiran2}*{Lemma 3.2}, any compact set in a complete subcritical Weinstein domain is Hamiltonian displaceable.

Let $U$ be an open neighbourhood of $\0$ in $W$.  Let $\delta > 0$ so that $F$ has no critical values in the interval $[- \delta, 0)$.  We then let $\widetilde W$ be the completion of the subcritical Weinstein domain $F^{-1}((-\infinity, -\delta])$.  We now observe that by following the flow of $Y$, we may identify  $F^{-1}(0) \setminus U$ with a compact set in the positive symplectization of $F^{-1}(\delta)$,  and thus with a compact set in $\widetilde W$.   The result now follows by the Lemma of Biran and Cieliebak.

\end{proof}

\section{The time dependent case}

Most of the elements of the proof for the autonomous case carry through to the time-independent case, 
but some new difficulties arise.  
Recall that by hypothesis, $(M, \omega, Y, F)$ is a Weinstein domain, with $Y$ a complete vector field.
Let $\lambda = i_Y \omega$ denote the corresponding primitive of the symplectic form.  
The Hamiltonian $H : S^1 \times M \to \R,\, (t, x) \mapsto H_t(x)$ 
defines a 1-periodic time-dependent Hamiltonian vector field $\Xt$.  
Furthermore, $dH_t( \0 ) = 0$ for all $t$, and $\Xt$ has a hyperbolic zero at $\0$ for all $t$.
By hypothesis, $\lambda(\Xt) - H_t \ge 0$, with equality only at $\0$.  
This condition implies that $S^1 \times (M \setminus \{ \0 \} )$ is a contact manifold with contact form
$\Lambda = \lambda - H_t dt$.  This, in turn, is a hypersurface of restricted contact type in $\R \times S^1 \times (M \setminus \{ \0 \})$.

As in the autonomous case, we construct a closed, embedded Lagrangian approximation 
to the unstable manifold.  The non-compactness of $M$ introduces some extra subtleties in carrying out the rest of the argument.  First of all, we will use the fact that $\lambda(\Xt) - H_t$ is proper to show that non-existence of homoclinic orbits contradicts Hamiltonian displaceability of the Lagrangian approximation to the unstable manifold.  Then, we will use the hypotheses that $H$ is positive outside a compact set, that $\lambda(X_{H_t}) - H_t$ and $H$ are both proper, and that $M$ is a Weinstein domain to construct a symplectic embedding in $\C \times M$, where the image has \textit{a priori} bounded displacement energy.

\begin{proof}[Proof of Theorem \ref{T:SelfCalibrated}]

The proof will again be by contradiction.  Assume there do not exist orbits homoclinic to $\0$.  

Let $\phi_t$ denote the flow of $\Xt$ on $M$.  With $t \in \R / \Z$, let $X = \die{t} + \Xt$.  
Denote the flow of $X$ by $\Sphi$.  Recall that $S^1 \times (M \setminus \{ \0 \})$ is a contact manifold 
with contact form $\Lambda = \lambda - H_t dt$.

By hypothesis, $\phi_1 : M \to M$ has a hyperbolic fixed point at $\0$.  We thus have the existence of the unstable manifold, an immersed $\R^n  \looparrowright  W^u_0 \subset M$.  We denote this immersion by $p_0$, with $p_0(0) = \0$.  As in the autonomous case, since $\0$ is a hyperbolic fixed point and since 
$D \phi_1(\0)$ is a linear symplectic map, the unstable manifold $W^u_0$ 
is an immersed Lagrangian submanifold.
Thus, there exists a function $g : W^u_0 \to \R$ with $g(0) = 0$, so that $p_0\pb\lambda = dg$.
Define now 
\begin{equation*}
\begin{aligned}
\tilde p : \R \times \R^n &\to S^1 \times M \\
(t, x) &\mapsto ( t \pmod 1, \phi_t \circ p_0(x) ).
\end{aligned}
\end{equation*}
Let 
\begin{equation*}
\widehat W^u = \R \times \R^n / \sim
\end{equation*}
where $(t, x) \sim (s, y)$ if $t \equiv s \pmod 1$ and $\phi_t( p_0(x)) = \phi_s(p_0(y))$.  
Then, $\tilde p$ descends to a smooth function $p : \widehat W^u \to S^1 \times M$.  
Let $W^u = p( \widehat W^u)$.  Observe that $X$ is tangent to $W^u$, so we may define 
$Z = p^{-1}\pf X|_{W^u}$ on $\widehat W^u$.  

Denote by $\FlowProj$ the map from $\R \times M$ to $S^1 \times M$ by following the flow of $X$.
\begin{equation}
\begin{aligned}
\FlowProj : \R \times M &\rightarrow S^1 \times M \\
(t, x) &\mapsto (t \mod 1, \phi_t(x) )
\end{aligned}
\end{equation}
Equivalently, this corresponds to $(t, x) \mapsto \Sphi_t( 0, x)$.

Let $\Action$ be the action of the trajectory through $x$, followed for time $t$.  
\begin{equation}
\begin{aligned}
\Action : \R \times M &\to \R \\
(t, x) &\mapsto \int_{0}^{t} 
		\phi_\tau \pb \left ( \lambda( X_{H_\tau} ) ) - H_\tau( \phi_\tau (x) \right )  d\tau.
	 \end{aligned}
\end{equation}
By hypothesis, $\Action$ is positive for $t > 0$ and $x \ne \0$.  
Observe that $\Action(t,x) = \int_0^t \Lambda(X)|_{\Sphi_\tau(0, x)} d\tau$. 
This corresponds to the right hand side of Equation \eqref{E:action} in the autonomous case.  
For fixed $t \in \R$, we will denote by $\Action_t : M \to \R$ the map $x \mapsto \Action(t, x)$.

We now have 
\begin{equation} \label{E:flowLambda}
\begin{split}
\phi_t \pb \lambda - \lambda &= \int_0^t \frac{d}{d \tau} \phi_\tau \pb \lambda d \tau \\
	&= \int_0^t  \phi_\tau \pb ( L_{X_{H_\tau}} \lambda ) d \tau \\
	&= \int_0^t d( \phi_\tau \pb ( \lambda(X_{H_\tau}) - H_\tau ) )d \tau \\
	&= d \Action_t(x).
\end{split}
\end{equation}

From this, it follows 
\begin{equation*}
\begin{split}
\FlowProj \pb \Lambda &= - \FlowProj \pb (H_t dt) + \FlowProj \pb \lambda \\
	&= - H_t( \phi_t(x) ) dt + \phi_t \pb \lambda + \lambda( \Xt ) dt \\
	&= \lambda + (\phi_t \pb \lambda - \lambda ) + 
			(- H_t( \phi_t(x) ) + \lambda( \Xt )|_{\phi_t(x)} ) dt \\
	&= \lambda + d \Action.
\end{split}
\end{equation*}
Hence, 
\begin{equation*}
\tilde p \pb \Lambda = d\left ( g + \Action( \cdot, p( \cdot ) )  \right ).
\end{equation*}

Considering \eqref{E:flowLambda} in the special case of $x \in W^u_0$, we also obtain 
\begin{equation} \label{E:g_action}
g( p_0^{-1}( \phi_1(x)) ) - g( p_0^{-1}(x) ) = \Action(1, x).
\end{equation}
We will now show that $\tilde G := g + \Action(\cdot, p(\cdot))$, 
defined on $\R \times W^u_0$ descends to a smooth function on $\widehat W^u$, 
which we denote by $G$.  Indeed, if $(t, x) \sim (s, y)$, 
we have that $s = t + k$, $k \in \Z$ and $x = \phi_k(y)$.  We then obtain 
\begin{equation*}
\begin{aligned}
\tilde G(s, y) &= g(y) + \Action(s, y) = g(y) + \Action(t + k, y) \\
	&= g(y) + \Action(k,  y) + \Action(t, x) \\
	&= g(x) + \Action(t, x) \\
	&= \tilde G(t, x).
\end{aligned}
\end{equation*}

We will now show that by the assumption of non-existence of homoclinic orbits, 
$G : \widehat W^u \to \R$ is proper.   The argument will be similar to the autonomous case.

First, by hypothesis, $\0$ is a hyperbolic fixed point of $\phi_1$.  The unstable manifold $W^u_0$ is 
invariant under $\phi_1$, so for each integer $k$, let 
$\hat \phi_k = p_0^{-1} \circ \phi_k \circ p_0 : \R^n \to \R^n$.
Then, on $\R^n$, there exist a norm $N$ and constants $r_2 > r_1 > 1$ so that for all $x \in \R^n$,
\begin{equation} \label{E:Hyperbolic}
r_1 N(x) \le N( \hat \phi_1 (x) ) \le r_2 N(x).
\end{equation}
Consider the annulus $Q = \{ x \, | \, 1 \le N(x) \le r_2 \}$.  Let $c = \inf \{ g(x) \, | \, x \in Q \}$.
Then, \[
\{ x \, | \, g(x) \le c \} \subset \{ x \, | \, N(x) \le r_2 \}
\]
and thus is compact.
Indeed, if for some $x$, $g(x) \le c$ but $N(x) > r_2$, we would have $y = \phi_{-1}(x)$ with
$g(y) < g(x) \le c$ and $r_1 N(y) \le N(x) \le r_2 N(y)$ and thus $1 \le N(y) \le r_2$.  

Let $0 < a < b$.  By the assumption of non-existence of homoclinic orbits, 
there is an open neighbourhood $U$ of $\0$ so that trajectories of the vector field $X$ 
starting on $W^u \setminus U$, will not enter $U$ for time $t \le b+1$.  
By shrinking $U$, we may assume that for all $x \in \R^n$ with $p(x) \in U$, $g( \phi_1(x) ) < \min \{ a, c \}$.  Let $\hat U = p_0^{-1}( U )$.
By hypothesis, $\lambda(X_{H_t}) - H_t$ is proper, and is positive on the complement of $U$.  
Thus, there exists a constant $c_U > 0$ so that $\lambda(X_{H_t}) - H_t \ge c_U$ on $M \setminus U$.
It follows then that for $x \in \R^n \setminus \hat U$, and for $k$ a positive integer,
\begin{align}
\lim_{k \to - \infinity} g( \hat \phi_k  (x) ) &= 0 \label{E:a}\\
g( \hat \phi_k (x) ) - g(x) &\ge c_U k. \label{E:b}
\end{align}
Suppose $x_n \in \R^n$ is a sequence with $g(x_n) \in [a, b]$ and with $g(x_n) \rightarrow g$.  
Then, by \eqref{E:a}, there exists a sequence of integers, $k_n$ so that 
\begin{equation*}
g( \hat \phi_{k_n} ( x_n ) ) \le c \quad \text{ and } \quad 
g( \hat \phi_{k_n + 1}  ( x_n ) ) > c.
\end{equation*}
Let $y_n = \hat \phi_{k_n} ( x_n )$.  
Now, by \eqref{E:b}, for each $n$, we obtain the following bounds :
\begin{align*}
\text{for $k_n \ge 1$ : } c - a &\ge g(  \hat \phi_{k_n} ( x_n ) ) - g(x_n) \ge c_U k_n \\
\text{for $k_n \le 1$ : } b - c &\ge g(  \hat \phi_{- k_n - 1} \circ \hat \phi_{k_n + 1} ( x_n )) 
		-	g( \hat \phi_{k_n + 1} ( x_n )) \ge - c_U ( k_n + 1).
\end{align*}
Thus, $| k_n| $ is bounded and, furthermore, the set of points with $g \le c$ is compact, so
there exists a subsequence $n_l$ so that $k_{n_l} \to k_*$ and $y_{n_l} \to y_*$.  
Then, $x_{n_l} \to \phi_{-k_*} ( y_* )$ is a convergent subsequence.
This shows that $g$ is proper.
Now, observe that $g \ge 0$ and $\Action( \cdot, p_0( \cdot ) ) \ge 0$ on $[0, 1] \times \R^n$.
Thus, $\tilde G = g + \Action( \cdot, p_0( \cdot ) )$ is proper on $[0, 1] \times \R^n$.  It follows that $G$ is proper on $\widehat W^u$.

Let $0 < a < b$, $S > 0$, and $U$ a neighbourhood of $\0$ in $M$, so that $S^1 \times (M\setminus U)$ is forward-flow invariant, and so that $G^{-1}(a) \cap U = \emptyset$.  
Then there exists an embedding :
\begin{equation*}
\begin{aligned}
\iext : [-S, 0] \times G^{-1}([a,b]) &\to [-S, 0] \times S^1 \times ( M \setminus U) \\
(s, x) &\mapsto (s, p(x) ).
\end{aligned}
\end{equation*}
Note that $[-S, 0] \times S^1 \times ( M \setminus U)$ is a symplectic manifold with symplectic form 
$d (\Exp{s} \Lambda)$.
The embedding satisfies 
\begin{equation*}
\iext \pb ( \Exp{s}\Lambda ) = \Exp{s}dG.
\end{equation*}

By Lemma \ref{L:embeds_in_C} (below), this may be symplectically embedded in $( \C \times M, \omega_0 + \omega)$.  Furthermore, with $\pi_1 : \C \times M \to \C$ denoting the projection, $\pi_1 \circ \iext$ has image contained in a compact set $K \subset \C$, with $K$ independent of $b$.  
There exists a Hamiltonian diffeomorphism $|| \Phi ||$ displacing $K$ from itself, with compact support and finite Hofer norm.  Let $B_\tau$ be a time dependent Hamiltonian on $\C$, generating $\Phi$, of oscillation no greater 
than $2 || \Phi ||$.  For any $b > 0$, the image of $[-S, 0] \times G^{-1}[a,b]$ lies in a compact region $\Omega \subset S^1 \times M$.  
Let $\beta$ be a smooth function with compact support, $\beta : \R \times S^1 \times M \to [0,1]$, with $\beta = 1$ on $\Omega$.  We now take the Hamiltonian on $\C \times M$ given by $\beta B_\tau$.  This now has compact support and will displace $G^{-1}[a,b]$.  Furthermore, its oscillation will be no greater than $2 || \tilde \Phi ||$.  

The result now follows by contradiction with Theorem \ref{T:Mohnke}.
\end{proof}

\begin{lemma}\label{L:embeds_in_C}
Let $(M, \omega, Y, F)$ be a complete Weinstein domain, with $\lambda = i_Y \omega$ the associated primitive of $\omega$.  Let $H_t : S^1 \times M \to \R$ be a one-periodic, time-dependent Hamiltonian so that 
$\lambda(\Xt) - H_t$ is proper and is positive away from $\0$, so $H_t$ is proper and positive for all $t$ outside of a compact region in $M$.  
Let $\omega_0$ denote the standard symplectic form on $\C$.  
Then, for any neighbourhood $U \subset M$ of $\0$, there exists a bijective symplectic immersion 
\[
i : \left ( (-\infinity, 0] \times S^1 \times (M \setminus U), d( \Exp{s} \Lambda ) \right ) 
	\to \left ( \C \times M , \omega_0 + \omega \right ).
\]
Furthermore, with $\pi_1 : \C \times M \to C$, denoting projection on to the first factor, $\pi_1 \circ i$ has bounded image.
\end{lemma}
\begin{proof}
The proof is by constructing the immersion.  This is done in three steps.
Let $\lambda = i_Y \omega$ be the primitive of $\omega$ on $M$.  

The first step is to smooth the contact form $\Lambda$ near $\0$, to obtain a contact form $\bar \Lambda$ on $S^1 \times M$.  
Indeed, by a construction in a Darboux chart, it is straightforward to find 
functions $P$ and $Q$ with support in $U$ so that with $\bar \lambda = \lambda + dP$
we have $\bar \lambda(X_Q) - Q \ge 0$ and $\bar \lambda(X_Q) - Q > 0$ at $\0$.
Let now $\bar H = H +Q$.  Define 
\begin{equation*} 
\bar \Lambda = \lambda' - (H + Q) d t.
\end{equation*}
Let $\bar X$ be the vector field $\die{t} + X_{\bar H_t}$.  We have $\bar \Lambda ( \bar X) > 0$ on $S^1 \times M$.  
By construction, $\bar \Lambda$ and $\Lambda$ agree outside of $U$, as do $\bar X$ and $X$.
Observe furthermore that with $\bar Y = Y - X_P$, we obtain a Liouville vector field so that 
$i_{\bar Y} \omega = \bar \lambda$.  Then, $dF[ \bar Y ] > 0$ outside of a compact region in $M$.

The second step is to show there   
exists a constant $c > 0$, a function $f : S^1 \times M \to \R$ bounded above, and a diffeomorphism 
\begin{equation*}
\Psi : S^1 \times M \to S^1 \times M \text{ with } 
\Psi \pb ( \bar \lambda + c dt) = \bar \Lambda = \bar \lambda - \bar H_t dt.
\end{equation*}
The proof of the second step is similar to that of Gray's theorem, which does not directly
apply due to the non-compactness of $S^1 \times M$.  
Observe that $\bar \lambda( \bar \Xt) - \bar H_t$ is positive and proper.
Thus, there exists $c > 0$ so that $c < \frac{1}{2}  \inf  \{ \bar \lambda( \bar \Xt )  - \bar H_t \, | \, (t, x) \in S^1 \times M \}$. 

Let $\bar H^\tau_t = \tau ( \bar H_t + c) - c$, and $\bar \Lambda^\tau = \bar \lambda - \bar H^\tau_t dt$.
Note that $\bar \Lambda^0 = \bar \lambda + c dt$ and 
	$\bar \Lambda^1 = \bar \Lambda = \bar \lambda - \bar H_t dt$.
Observe also $X_{\bar H^\tau_t} = \tau X_{\bar H_t}$.  

Define the following function on $S^1 \times M$ :
\begin{equation*}
h^\tau(t, x) =  - \frac{\bar H_t(x) + c}{\tau \left ( \alpha(X_{\bar H_t}) - \bar H_t -c \right ) + c}.
\end{equation*}
Define also the following vector field on $S^1 \times M$ : 
\begin{equation*}
Z^\tau = h^\tau \bar Y.
\end{equation*}
Let $\Phi^\tau$ be the flow of $Z^\tau$ on $S^1 \times M$.  
By hypothesis, $\bar H_t$ is positive outside a compact region in $M$, 
so $h_\tau$ will be negative outside a compact set.  
Since $dF[ \bar Y ] > 0$ outside a compact set, 
it follows that $L_{Z_\tau} F$ is negative outside a compact set.  
By hypothesis, $F$ is bounded below, so the flow of $Z_\tau$ exists for all $\tau \in [0, 1]$.

By the construction of $h^\tau$, and Cartan's formula, we obtain :
\begin{equation*}
\begin{split}
\frac{d}{d\tau} \left( \Phi^\tau \pb \bar \Lambda^\tau \right ) 
		&= \Phi^\tau \pb \left ( \omega( Z^\tau, \cdot ) 
				- d\bar H^\tau_t( Z^\tau ) dt + \frac{d}{d\tau} \bar \Lambda^\tau \right ) \\
		&= \Phi^\tau \pb \left ( h^\tau \cdot \left 
				( \bar \lambda 
				- \left ( 
					\bar \lambda( X_{\bar H^\tau_t} ) + \frac{\bar H_t + c}{h^\tau} 
					\right ) dt 
				) \right ) 
				\right )\\
		&= \Phi^\tau \pb ( h^\tau \bar \Lambda^\tau ).
\end{split}
\end{equation*}
By integrating this equation, we obtain :
\begin{equation*}
\Phi^1 {\pb} \Lambda^1 = \Exp{f} \Lambda^0 ,\text{ with } 
 f = \int_0^1 h_\tau \circ \Phi^\tau d \tau .
\end{equation*}
Since $h_\tau$ is negative outside a compact set, it is bounded from above.  Thus, $f$ is also bounded from above.
This completes the proof of the second step, with $\Psi = \Phi^1$.

The third and final step is to show that 
$(S^1 \times M, c dt + \bar \lambda)$ can be realized as a hypersurface of restricted contact type in 
$\C \times M$.  
Let $\psi_t$ denote the flow of $\bar Y$ on $M$.  
Recall that our convention is $S^1 = \R / \Z$.  
Consider the following embedding :
\begin{align*}
i' : S^1 \times M & \to \C \times M \\
(t, x) &\mapsto \left ( \frac{c}{\pi} \Exp{\frac{1}{2} f(t, x) + 2 \pi i t}, \psi_{f(t,x)} (x) \right ).
\end{align*}
On $\C$, with coordinates $x + iy$, we have the standard symplectic form $\omega_0 = d \lambda_0$, where $\lambda_0 = 1/2( x dy - y dx)$.  
A computation shows that $i' \pb (\lambda_0 + \lambda) = \Exp{f} \left ( \lambda + c dt \right )$. 
This embedding extends to a symplectic embedding :
\begin{align*}
\iext' : (-\infinity, 0] \times S^1 \times M & \to \C \times M \\
(s, t, x) &\mapsto \left ( \frac{c}{\pi} \Exp{\frac{1}{2} s + \frac{1}{2} f(t, x) + 2 \pi i t}, 
				\psi_{s + f(t,x)} (x) \right ).
\end{align*}

By hypothesis, $f$ is bounded above.  Hence $\Exp{ \frac{1}{2} ( f + s) }$ is bounded for $s \le 0$, and thus the image of $\pi \circ i' $ is bounded.  

The desired immersion is then given by $i = \iext' \circ \Psi$.
\end{proof}



\def\cprime{$'$}


\end{document}